\documentclass[11pt,english,oneside]{amsart}

\usepackage[breaklinks,colorlinks,linkcolor=blue,citecolor=blue,urlcolor=blue]{hyperref}
\usepackage[T1]{fontenc}
\allowdisplaybreaks[4]
\usepackage[latin9]{inputenc}
\usepackage{geometry}
\usepackage{indentfirst}
\usepackage{amssymb}
\usepackage{color}
\usepackage{float}
\usepackage{algorithmic}
\usepackage{algorithm}
\usepackage[caption=false]{subfig}
\usepackage{xcolor}

\usepackage{graphicx}
\usepackage{graphics}
\usepackage{epsfig}
\usepackage{graphicx,color}
\usepackage{amsmath, amssymb}

\makeatletter
\numberwithin{equation}{section} 
\numberwithin{figure}{section} 
  \@ifundefined{theoremstyle}{\usepackage{amsthm}}{}
  \theoremstyle{plain}
  
  \theoremstyle{plain}
  
  \theoremstyle{plain}
  
  \theoremstyle{remark}
  
  \theoremstyle{remark}
  
  \theoremstyle{plain}
  
\usepackage{geometry}

\def\e{\hbox{\rm e}}

\def\exp{\hbox{\rm exp}}
\smallskip
\def\<{{\langle }}
\def\>{{\rangle }}

\makeatother

\usepackage{babel}

\def\e{\hbox{\rm e}}

\def\exp{\hbox{\rm exp}}
\smallskip
\def\<{{\langle }}
\def\>{{\rangle }}

\makeatother

\usepackage[mathscr]{eucal}
\usepackage{amssymb}
\usepackage{latexsym}
\usepackage{amsthm}
\theoremstyle{plain}
\newtheorem{theorem}{Theorem}[section]

\newtheorem{remark}{Remark}[section]

\newtheorem{definition}{Definition}[section]
\newtheorem{question}{Question}[section]

\makeatletter
\@addtoreset{equation}{section}

\usepackage{color}

\usepackage{graphicx,color}

\title[ Colding-Minicozzi entropies of some self-shrinkers]
{Colding-Minicozzi entropies of some self-shrinkers}
\author{Qilun Luo, Guoxin Wei and Fu-An Zhang}
\address{Qilun Luo \\  School of Mathematical Sciences, South China Normal University,
510631, Guangzhou,  China, qilunluo92@gmail.com}
\address{Guoxin Wei \\  School of Mathematical Sciences, South China Normal University,
510631, Guangzhou,  China, weiguoxin@tsinghua.org.cn}
\address{Fu-An Zhang \\  School of Mathematical Sciences, South China Normal University,
510631, Guangzhou,  China, cvgmt@163.com}

\begin{document}
\maketitle

\begin{abstract}
In this note, we numerically estimate Colding-Minicozzi entropies of some self-shrinkers and get that Colding-Minicozzi entropies of $n$-dimensional Angenent torus are decreasing about dimension $n$ ($2\leq n\leq 5*10^7$), which partially answer the questions of Berchenko-Kogan \cite{BK}.
\end{abstract}


\section{Introduction}

\begin{definition}
A hypersurface $X: M^n\rightarrow \mathbb{R}^{n+1}$ is called a self-shrinker  of mean curvature flow if it satisfies the equation
\begin{equation}\label{eq:2022-12-1}
H+\langle X, N \rangle=0,
\end{equation}
where $X$ is the position vector, $H$ is the mean curvature, and $N$ is the unit normal vector.
\end{definition}

Self-shrinkers are key to understand the singularities that develop under mean curvature flow, we concentrate on the so-called "Type I" singularity, studied by Huisken in \cite{H2}. Self-shrinkers plays an important role for study on singularities of the mean curvature flow. Mean curvature flow is the important geometric flow.

There are some trivial self-shrinkers. Examples of self-shrinkers in $\mathbb{R}^{n+1}$ include the $n$-dimensional Euclidean space $\mathbb{R}^{n}$, the $n$-dimensional sphere $S^n(\sqrt{n})$ of radius $\sqrt{n}$ centered at the origin, the $n$-dimensional cylinder  $S^m(\sqrt{m})\times \mathbb{R}^{n-m}$  for $1\leq m\leq n-1$.

In 1989,  Angenent  \cite{A} constructed compact embedded self-shrinker
$X:S^1\times S^{n-1}\to\mathbb{R}^{n+1}$ which is the first nontrivial self-shrinkers, a torus. Later, the torus is called Angenent torus. Angenent  \cite{A} also provided the numerical evidence for the existence of an immersed sphere self-shrinker. Until now, there are few geometric and topological properties about Angenent torus, such as, the uniqueness of the shrinking doughnuts is still open in all dimensions $n \geq 2$. For self-shrinkers, there are a well-known conjecture asserts that the round sphere  should be the only embedded self-shrinker which is diffeomorphic to a sphere. Brendle \cite{B} proved the above conjecture for $2$-dimension self-shrinker. For the higher dimensional  self-shrinker, the conjecture is still open.

Many interesting  examples of complete self-shrinker were found recently.
For instance,  complete embedded self-shrinkers with higher genus in $\mathbb{R}^{3}$ were constructed by
Kapouleas, Kleene and M{\o}ller \cite{KKM} (also see \cite{KM}, \cite{M} ) and Nguyen \cite{N1}-\cite{N3}.
Drugan \cite{D} constructed an immersed self-shrinker, which is a topological sphere
$X:S^2\to\mathbb{R}^{3}$. In \cite{DK}, Drugan and Kleene constructed many complete immersed self-shrinkers with rotational symmetry for each of following topological types: the sphere, the plane, the cylinder and the torus. In \cite{CW1}, Cheng and Wei constructed some special immersed rotational compact $2$-dimensional self-shrinkers whose profile curves do not intersect symmetry axis perpendicularly; see Figure 1.1.

\begin{figure}[ht]
\centerline{
\ \ \ \ \ \ \ \  \includegraphics[height=2cm]{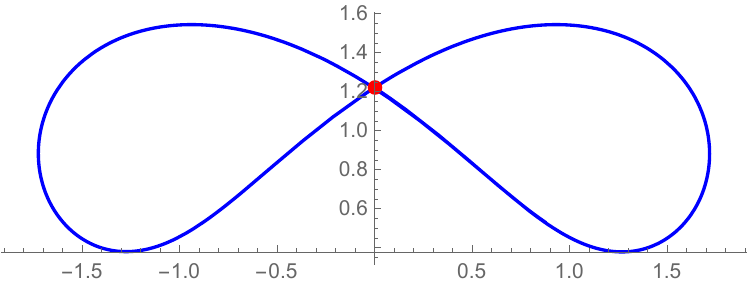}
\includegraphics[height=4.4cm]{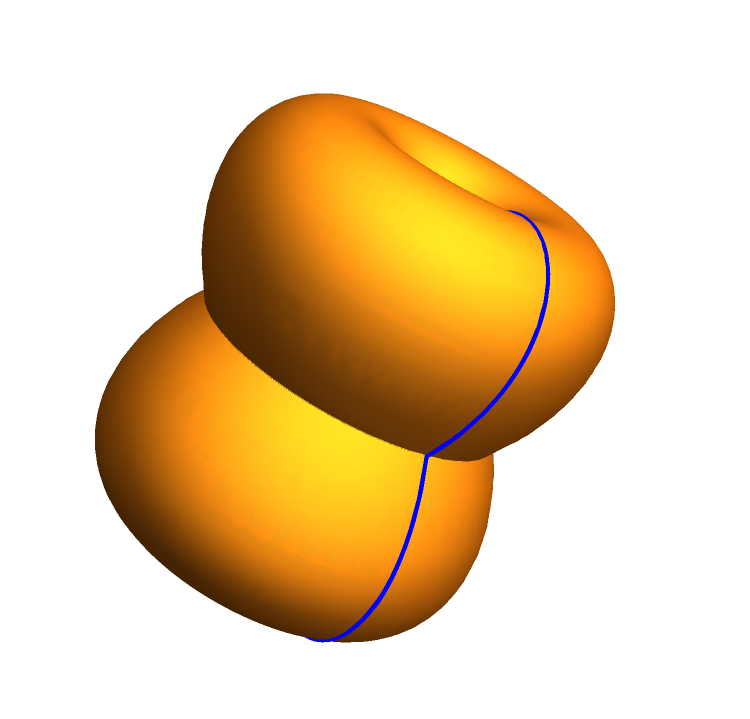}
\includegraphics[height=4.4cm]{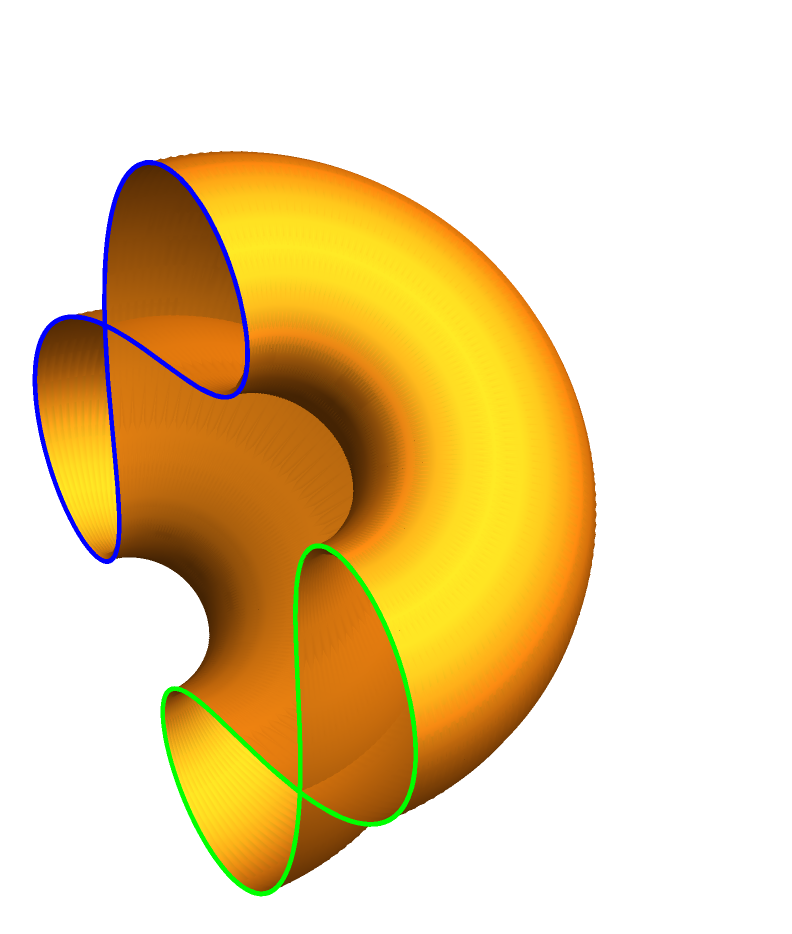}
}

\caption{ The graph of profile curve of compact self-shrinker,
self-shrinker and half of self-shrinker, here $n=2$. }

\label{thegraphs}
\end{figure}

 Huisken studied the singularities of mean curvature flow using a monotonicity result for a Gaussian-weighted area called the $F$-functional. The $F$-functional is defined by

\begin{definition}
For $X: M^n\rightarrow \mathbb{R}^{n+1}$, the $F$-functional $F(M)$ is
\begin{equation}
F_{X_0,t_0}(M)=(2 \pi t_0)^{-n/2}\int_M e^{-|X-X_0|^2/(2t_0)}d\mu.
\end{equation}
\end{definition}
\noindent Here $(2 \pi t)^{-n/2}e^{-|X|^2/(2t)}$ is the heat kernel of $\mathbb{R}^n$. Colding and Minicozzi \cite{CM} introduced the following functional which is called the entropy of $X: M^n\rightarrow \mathbb{R}^{n+1}$
on the space of $n$-dimensional submanifolds  $X: M^n\rightarrow \mathbb{R}^{n+1}$.

\begin{definition}
For $X: M^n\rightarrow \mathbb{R}^{n+1}$, the entropy $\lambda[M]$ is
\begin{equation}
\lambda[M]=\sup\limits_{X_0,t_0}  F_{X_0,t_0}(M).
\end{equation}
\end{definition}

The entropy is a natural geometric quantity which measures the complexity of a hypersurface in $\mathbb{R}^{n+1}$. Colding and Minicozzi observed that, by Huisken's monotonicity formula \cite{H2}, the entropy of $X: M^n\rightarrow \mathbb{R}^{n+1}$ is non-increasing along the mean curvature flows and plays a significant role in analyzing the dynamics of this flow. It is very interesting to study the entropies of hypersurfaces, especially, the entropies of self-shrinkers. For research on self-shrinkers, see \cite{CL, CW, CM, DX,  H2}.

For a self-shrinker, entropy remains constant during mean curvature flow, that is, $\lambda[M]=F_{o,1}(M)$; for any other hypersurface, the entropy is monotonically decreasing during mean curvature flow. In \cite{S}, Stone computed the entropies of $S^m(\sqrt{m})\times \mathbb{R}^{n-m}$ and obtained that the values of the entropies of $S^m(\sqrt{m})\times \mathbb{R}^{n-m}$ are all distinct, and bigger than $1$ for $m>0$. Indeed, the entropies form a strictly decreasing sequence in $m$, with the entropies tending to $\sqrt{2}$ as $m\rightarrow\infty$. In \cite{CWZ}, Cheng-Wei-Zeng obtained that if $M$ is a compact minimal rotational hypersurface in $S^{n+1}(1)$, then the entropies $\lambda[C(M)]$ of the cones $C(M)$ over $M$ in $\mathbb{R}^{n+2}$ satisfies either $\lambda[C(M)]=1$, or $\lambda[C(M)]=\frac{2\pi \sigma_{n-1}\sqrt{a_0}}{\sigma_n}$, or
$\lambda[C(M)]>\frac{4(\pi-1) \sigma_{n-1}\sqrt{a_0}}{\sigma_n}$, here $a_0=\frac{(n-1)^{n-1}}{n^n}$, $\sigma_n$ denotes the $n$-area of $S^n(1)$. In \cite{CS}, Chu and Sun got that a smooth embedded self-shrinker in $\mathbb{R}^{3}$ with the fourth lowest entropy is not rotationally symmetric.

In 2013, Colding-Ilmanen-Minicozzi-White \cite{CIMW} showed that within the class of compact smooth self-shrinkers of the mean curvature flow in $\mathbb{R}^{n+1}$, the entropy is uniquely minimized at the round sphere. They conjectured that, for $2\leq n\leq 6$, the round sphere minimizes the entropy among all compact hypersurfaces. Using an appropriate weak mean curvature flow, Bernstein and Wang \cite{BW} proved their conjecture. The $n=2$ and non-toric case was also proved independently by Ketover and Zhou \cite{KZ} using a min-max theory for the Gaussian surface area. Zhu \cite{Z} proved that the round sphere minimizes the entropy among all compact embedded  hypersurfaces for all dimensions $n$ by using an extension of Colding-Minicozzi's classification of entropy-stable self-shrinkers to the singular setting.


Mramor \cite{Mr} calculated that the rotationally symmetric shrinking donuts of Drugan and Nguyen \cite{DN} have entropy less than $2$. Ma, Muhammad and M{\o}ller \cite{MMM} got that the entropy of complete embedded rotationally symmetric self-shrinker is bound by a explicit number $E_n$, where $2.02780<E_n<2.24759$. In 2019, Berchenko-Kogan \cite{BK} numerically estimated the entropy of $2$-dimensional Angenent torus using the discrete Euler-Lagrange equations and obtained that the entropy of $2$-dimensional Angenent torus is approximately 1.85122. In the same paper, he also proposed the following question.

\begin{question}
What are the entropies of $n$-dimensional $(n>2)$ Angenent torus? What value does the entropy approach as $n$ becomes large?
\end{question}

In this note, we numerically estimate the Colding-Minicozzi entropy of $n$-dimensional ($2\leq n\leq 5*10^7$) Angenent torus by programming and partially answer the questions of Berchenko-Kogan.

\begin{theorem}

For $2\leq n\leq 5*10^7$, the entropies of $n$-dimensional Angenent torus  are decreasing about dimension $n$, the perimeters of profile curves of $n$-dimensional Angenent torus are also decreasing about dimension $n$. Moreover, the entropy of $2$-dimensional Angenent torus is approximately $1.85121667$, the perimeter of profile curve of $2$-dimensional Angenent torus is approximately $5.30925757$.

\end{theorem}

\begin{remark}
Using our method, we can  get the entropy of $2$-dimensional round sphere is approximately $1.471517764$, whereas the true value is $\frac{4}{\e}\approx 1.471517765$.
\end{remark}

\begin{remark}
Berchenko-Kogan \cite{BK} conjectured that the values of entropies about the $2$-dimensional Angenent torus appear to lie on an exponential curve converging to $1.8512167$. In fact, we get the entropy of $2$-dimensional Angenent torus is approximately $1.85121667$.
\end{remark}

%

In addition, we also numerically estimate Colding-Minicozzi entropy of self-shrinkers constructed by Mcgrath \cite{M}, self-shrinkers constructed by Cheng-Wei \cite{CW1}.

\begin{theorem}

For $2\leq m\leq 5*10^7$, the entropies of $(2m-1)$-dimensional self-shrinkers constructed by Mcgrath \cite{M} are decreasing about dimension. Moreover, the entropy of $3$-dimensional this self-shrinker is approximately $2.46576946$ $(>E_3)$. For $2\leq n\leq 5*10^7$, the entropies of $n$-dimensional self-shrinkers constructed in \cite{CW1} are also decreasing about dimension $n$. Moreover, the entropy of $2$-dimensional this self-shrinker is approximately $2.88472911$.
\end{theorem}

\section{Preliminaries}

In \cite{A}, Angenent used a "shooting method" to construct "Angenent's torus" solution to mean curvature flow.  In fact,
Let $\gamma(s)=(x(s), r(s))$, $s\in (a,b)$ be a curve with $r>0$ in
the upper half plane $\mathbb{H}=\{x+ir| \  r>0, \  x\in \mathbb{R}, \ i=\sqrt{-1}\}$,
where $s$ is arc length parameter of $\gamma(s)$.
We consider  a rotational hypersurface $X: (a,b)\times S^{n-1}(1) \hookrightarrow \mathbb{R}^{n+1}$
in $\mathbb{R}^{n+1}$ defined by
\begin{equation}\label{eq:10-30-1}
X:  (a,b)\times S^{n-1}(1)\hookrightarrow \mathbb{R}^{n+1}, \ \ X(s, \alpha)=(x(s), r(s)\alpha) \in \mathbb{R}^{n+1}
\end{equation}
where  $S^{n-1}(1)$ is the ($n-1$)-dimensional unit sphere (cf. \cite{DD}).

\noindent
By a direct calculation, one has the unit normal vector
\begin{equation}\label{eq:10-30-2}
N=(-r^{\prime}, x^{\prime}\alpha)
\end{equation}
and the mean curvature
\begin{equation}\label{eq:10-30-3}
H=-x^{\prime\prime}r^{\prime}+x^{\prime}r^{\prime\prime}-\dfrac{n-1}{r}x^{\prime}.
\end{equation}
Therefore, we know from \eqref{eq:10-30-1} and \eqref{eq:10-30-2}

\begin{equation}\label{eq:10-30-4}
\langle X, N\rangle=-xr^{\prime}+rx^{\prime}.
\end{equation}
Hence,  $X:  (a,b)\times S^{n-1}(1)\hookrightarrow \mathbb{R}^{n+1}$ is a self-shrinker in $\mathbb{R}^{n+1}$,
if and only if,  from \eqref{eq:2022-12-1}, \eqref{eq:10-30-3} and \eqref{eq:10-30-4},
\begin{equation}\label{eq:9-16-11}
-x^{\prime\prime}r^{\prime}+x^{\prime}r^{\prime\prime}-\dfrac{n-1}{r}x^{\prime}-xr^{\prime}+rx^{\prime}=0.
\end{equation}
Since $s$ is arc length parameter of the profile curve $\gamma(s)=(x(s), r(s))$, we have
\begin{equation}\label{eq:9-4-11}
(x^{\prime})^2+(r^{\prime})^2=1,
\end{equation}
Thus, it follows that
\begin{equation}\label{eq:9-16-12}
x^{\prime}x^{\prime\prime}+r^{\prime}r^{\prime\prime}=0.
\end{equation}
The signed curvature $\kappa(s)$ of the profile curve $\gamma(s)=(x(s), r(s))$ is given by
\begin{equation}\label{eq:12-4-2}
\kappa(s)=-\dfrac{x^{\prime\prime}}{r^{\prime}}
\end{equation}
and it is known that the integral of the signed curvature $\kappa(s)$ measures the total rotation of the tangent vector of
 $\gamma(s)$. From \eqref{eq:9-16-11} and \eqref{eq:9-16-12}, one has
\begin{equation}\label{eq:9-4-2}
x^{\prime\prime}=-r^{\prime}\bigl[xr^{\prime}+\bigl(\dfrac{n-1}{r}-r\bigl)x^{\prime}\bigl].
\end{equation}
Hence, we have

\begin{equation}\label{eq:9-4-3}
  {\begin{cases}
     (x^{\prime})^2+(r^{\prime})^2=1\\[2mm]

     -\dfrac{x^{\prime\prime}}{r^{\prime}}=xr^{\prime}+\bigl(\dfrac{n-1}{r}-r\bigl)x^{\prime}.
     \end{cases}}
  \end{equation}

  Let $x^{\prime}(s)=\cos\theta(s)$, then $r^{\prime}(s)=\sin\theta(s)$ and
  \begin{equation}\label{eq:d-theta}
  \theta^{\prime}(s)=x(s) \sin\theta(s)+((n-1)/r(s)-r(s)) \cos\theta(s).
  \end{equation}

Hence, the entropies of compact rotational self-shrinkers $X: M^n\rightarrow\mathbb{R}^{n+1}$, including Angenent torus, can be computed by
\begin{equation}
\lambda [M]=(2 \pi)^{-n/2}\int_M \exp^{-|X|^2/2} d\mu=(2 \pi)^{-n/2}\sigma_{n-1}\int_L r^{n-1}\exp^{-(x^2+r^2)/2} ds,
\end{equation}
where $\sigma_{n-1}$ denotes the area of $(n-1)$ dimensional unit sphere, $L$ is the perimeter of profile curve of rotational self-shrinkers.

We next compute the perimeter of profile curve of rotational self-shrinkers and obtain the entropies of compact rotational self-shrinkers $X: M^n\rightarrow\mathbb{R}^{n+1}$ by Algorithm \ref{alg:x0}.

\renewcommand{\algorithmicrequire}{\textbf{ Input:}}
\renewcommand{\algorithmicensure}{\textbf{Output:}}
\begin{algorithm}[ht]
\caption{ Compute the perimeters and entropies for the case: Angenent torus}
\label{alg:x0}
\begin{algorithmic}[1]
\REQUIRE the dimension $n$
\STATE \textbf{Initialize} $r_L = 10^{-5}, r_R = \sqrt{n-1}, L_{max} = 6, tspan = [0,L_{max}]$.
\STATE Set $opts = \mathrm{odeset}(``Events", @eventFun, ``RelTol", 10^{-10})$.
\WHILE{$|r_L-r_R|>10^{-10}$}
\STATE $r_0 = \frac{r_L+r_R}{2}$
\STATE $\mathbf{y_0} = [0;~ r_0;~ 0; ~0]$
\STATE Solve \eqref{eq:d-theta} by the high-order Runge-Kutta solver \emph{ode89}, i.e.,
\[
[\textbf{t}, \textbf{x}] = \emph{ode89}(@\mathrm{func}, tspan, \mathbf{y_0}, opts),
\]
where
$
\mathrm{func}(\textbf{t}, \textbf{y}) =  \left[ \cos(y_3);\sin(y_3); y_1\sin(y_3)+(\frac{(n-1)}{y_2}-y_2)\cos(y_3); \frac{2^{1-\frac{n}{2}}y_2^{n-1}e^{-\frac{y_1^2+y_2^2}{2}}}{\Gamma(\frac{n}{2})}\right].
$
\IF{$any(\mathbf{x_{:,2}})<r_0$}
\STATE $r_R = r_0$
\ELSE
\STATE $r_L = r_0$
\ENDIF
\ENDWHILE
\STATE Set the perimeter $L = \textbf{t}_{end}$
\STATE Set the entropy $\lambda = \textbf{x}_{end,4}$
\ENSURE the perimeter $L$ and the entropy $\lambda$
\end{algorithmic}
\end{algorithm}

\begin{table}[H]
    \centering
    \caption{Numerical results for the case: Angenent torus}
    \label{tab:Angenent-torus-new}
    \resizebox{0.65\linewidth}{!}{
    \begin{tabular}{r|c|c}
    \noalign{\hrule height 1pt}
  Dimension & Perimeter of profile curve  & Entropy  \\
  \noalign{\hrule height 1pt}
2		   & 5.30925757 &	1.85121667  \\ \hline
3		   & 5.27363687 &	1.80277855  \\ \hline
4		   & 5.26292303 &	1.78388334  \\ \hline
5		   & 5.25776364 &	1.77399119  \\ \hline
10		   & 5.24944377 &	1.75703279  \\ \hline
30		   & 5.24499759 &	1.74754529  \\ \hline
60		   & 5.24399424 &	1.74537297  \\ \hline
100		   & 5.24360376 &	1.74452494  \\ \hline
300		   & 5.24321928 &	1.74368858  \\ \hline
500		   & 5.24314309 &	1.74352269  \\ \hline
1000	   & 5.24308610 &	1.74339858  \\ \hline
3000	   & 5.24304818 &	1.74331598  \\ \hline
5000	   & 5.24304060 &	1.74329947  \\ \hline
10000	   & 5.24303492 &	1.74328710  \\ \hline
30000	   & 5.24303113 &	1.74327885  \\ \hline
50000	   & 5.24303038 &	1.74327720  \\ \hline
100000	   & 5.24302981 &	1.74327596  \\ \hline
300000	   & 5.24302943 &	1.74327514  \\ \hline
500000	   & 5.24302935 &	1.74327497  \\ \hline
1000000	   & 5.24302930 &	1.74327485  \\ \hline
3000000	   & 5.24302926 &	1.74327477  \\ \hline
5000000	   & 5.24302925 &	1.74327475  \\ \hline
10000000   & 5.24302925 &	1.74327472  \\ \hline
30000000   & 5.24302924 &	1.74327471  \\ \hline
50000000   & 5.24302924 &	1.74327469  \\ 
\noalign{\hrule height 1pt}
\end{tabular}
}
\end{table}


\section{Simulations}

\begin{figure}[H]
\includegraphics[width = 0.9\textwidth]{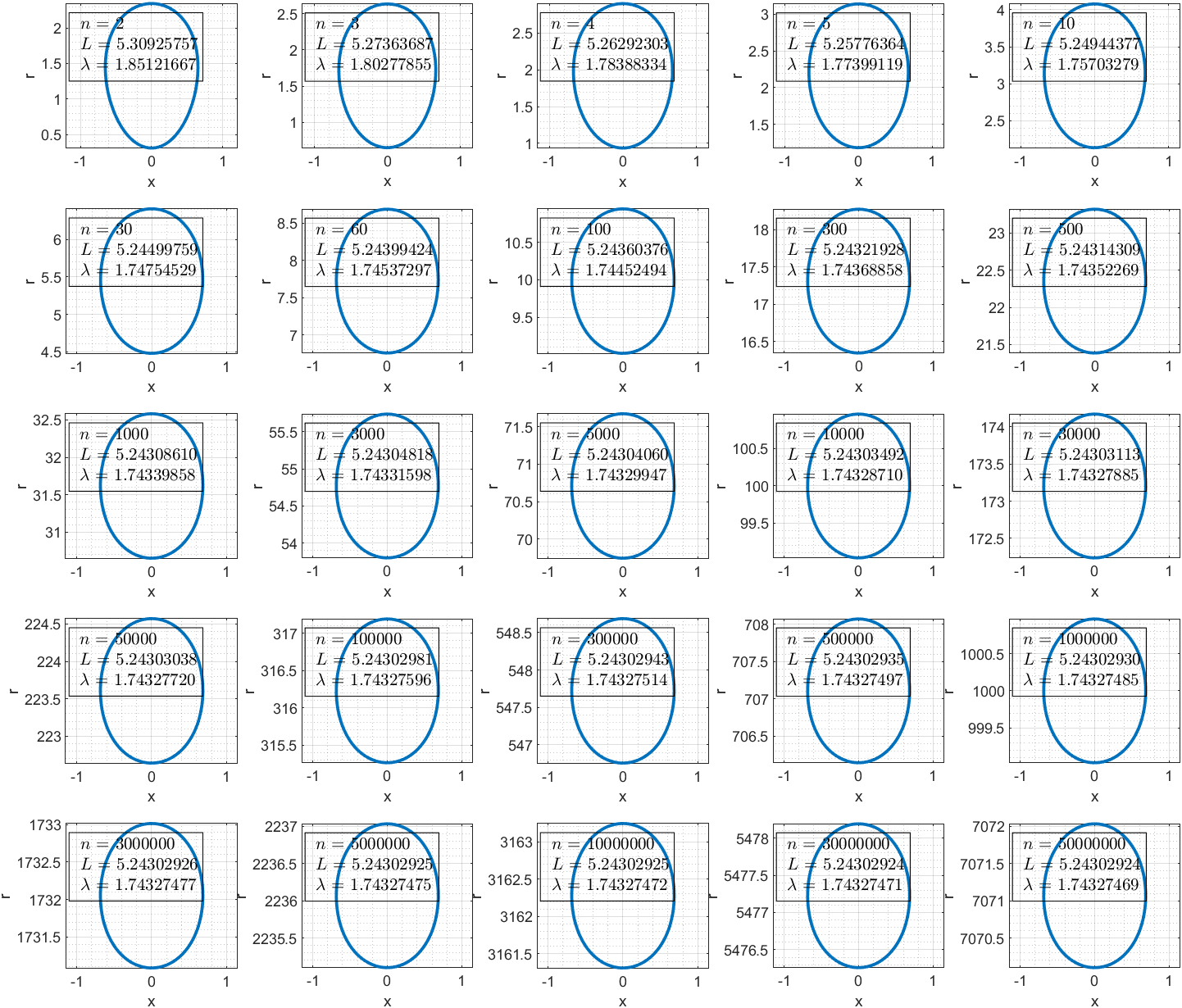}
\caption{ The graph of profile curve of Angenent torus for the case: Angenent torus. }
\label{fig:Angenent_torus_new}
\end{figure}



For each $n\geq 2$, inspired by constructions of Hsiang, Mcgrath \cite{M} constructed a new closed embedded self-shrinker in $\mathbb{R}^{2m}$. These self-shrinkers are diffeomorphic to $S^{m-1}\times S^{m-1}\times S^1$ and are $SO(m)\times SO(m)$ invariant. The entropies of these self-shrinkers $X: M^{2m-1}\rightarrow\mathbb{R}^{2m}$ can be computed by
\begin{equation}
\lambda [M]=(2 \pi)^{-\frac{2m-1}{2}}\int_M \exp^{-\frac{|X|^2}{2}} d\mu=\dfrac{\pi^{\frac{1}{2}}}{  2^{m-\frac{5}{2}}(\Gamma(\frac{m}{2}))^{2}}\int_L x^{m-1}r^{m-1}\exp^{-\frac{x^2+r^2}{2}} ds,
\end{equation}
where $\Gamma$ denotes the Gamma function, $L$ is the perimeter of profile curve of self-shrinkers constructed by Mcgrath.

 In this note,  we also numerically estimate the entropy of $(2m-1)$-dimensional $(2 \leq m \leq 5*10^7)$ these self-shrinkers by Algorithm \ref{alg:Mcgrath}.

\renewcommand{\algorithmicrequire}{\textbf{ Input:}}
\renewcommand{\algorithmicensure}{\textbf{Output:}}
\begin{algorithm}[ht]
\caption{ Compute the perimeters and entropies for the case: self-shrinkers constructed by Mcgrath}
\label{alg:Mcgrath}
\begin{algorithmic}[1]
\REQUIRE $m$
\STATE \textbf{Initialize} $r_L = 10^{-5}, r_R = \sqrt{m-1}, L_{max} = 6, tspan = [0,L_{max}]$.
\STATE Set $opts = \mathrm{odeset}(``Events", @eventFun, ``RelTol", 10^{-10})$.
\WHILE{$|r_L-r_R|>10^{-10}$}
\STATE $r_0 = \frac{r_L+r_R}{2}$
\STATE $\mathbf{y_0} = [r_0;~ r_0;~ \frac{7}{4}\pi; 0]$
\STATE Solve the ODEs by the high-order Runge-Kutta solver \textit{ode89}, i.e.,
\[
[\textbf{t}, \textbf{x}] = \emph{ode89}(@\mathrm{func}, tspan, \mathbf{y_0}),
\]
\[
\begin{aligned}
    \mathrm{func}(\textbf{t}, \textbf{y}) =  &\left[  \cos(y_3);\sin(y_3);(y_1-\frac{m-1}{y_1})\sin(y_3)+(\frac{(m-1)}{y_2}-y_2)\cos(y_3); \right.\\
     &~\left.\frac{4\pi^{m}(y_1y_2)^{m-1}e^{-\frac{y_1^2+y_2^2}{2}}}{(2\pi)^{\frac{2m-1}{2}}\Gamma(\frac{m}{2})^2}\right].
\end{aligned}
\]
\IF{$any(x(:,1)+x(:,2))<2r_0$}
\STATE $r_R = r_0$
\ELSE
\STATE $r_L = r_0$
\ENDIF
\ENDWHILE
\STATE Set the perimeter $L = \textbf{t}_{end}$
\STATE Set the entropy $\lambda = \textbf{x}_{end,4}$
\ENSURE the perimeter $L$ and the entropy $\lambda$
\end{algorithmic}
\end{algorithm}


\begin{table}[ht]
    \centering
    \caption{Numerical results for the case: self-shrinkers constructed by Mcgrath}
    \label{tab:Mcgrath-new}
    \resizebox{0.6\linewidth}{!}{
    \begin{tabular}{r|c|c}
    \noalign{\hrule height 1pt}
  m & Perimeter of profile curve  & Entropy  \\
  \noalign{\hrule height 1pt}
2          &    4.43826945 	& 2.46576946    \\ \hline
3          &    4.44243932 	& 2.31878674    \\ \hline
4          &    4.44299929 	& 2.26407546    \\ \hline
5          &    4.44312496 	& 2.23590016    \\ \hline
10         &    4.44310546 	& 2.18824431    \\ \hline
30         &    4.44297176 	& 2.16189319    \\ \hline
60         &    4.44292883 	& 2.15588727    \\ \hline
100        &    4.44291082 	& 2.15354529    \\ \hline
300        &    4.44289235 	& 2.15123696    \\ \hline
500        &    4.44288861 	& 2.15077928    \\ \hline
1000       &    4.44288579 	& 2.15043689    \\ \hline
3000       &    4.44288391 	& 2.15020903    \\ \hline
5000       &    4.44288353 	& 2.15016350    \\ \hline
10000      &    4.44288324 	& 2.15012936    \\ \hline  
30000      &    4.44288306 	& 2.15010661    \\ \hline  
50000      &    4.44288302 	& 2.15010206    \\ \hline  
100000     &    4.44288299 	& 2.15009864    \\ \hline
300000     &    4.44288297 	& 2.15009637    \\ \hline
500000     &    4.44288297 	& 2.15009592    \\ \hline
1000000    &    4.44288296 	& 2.15009558    \\ \hline
3000000    &    4.44288296 	& 2.15009536    \\ \hline
5000000    &    4.44288296 	& 2.15009529    \\ \hline
10000000   &    4.44288296 	& 2.15009523    \\ \hline  
30000000   &    4.44288296 	& 2.15009518    \\ \hline  
50000000   &    4.44288296 	& 2.15009512    \\ 
\noalign{\hrule height 1pt}
\end{tabular}
}
\end{table}


\begin{figure}[H]
\includegraphics[width = 0.8\textwidth]{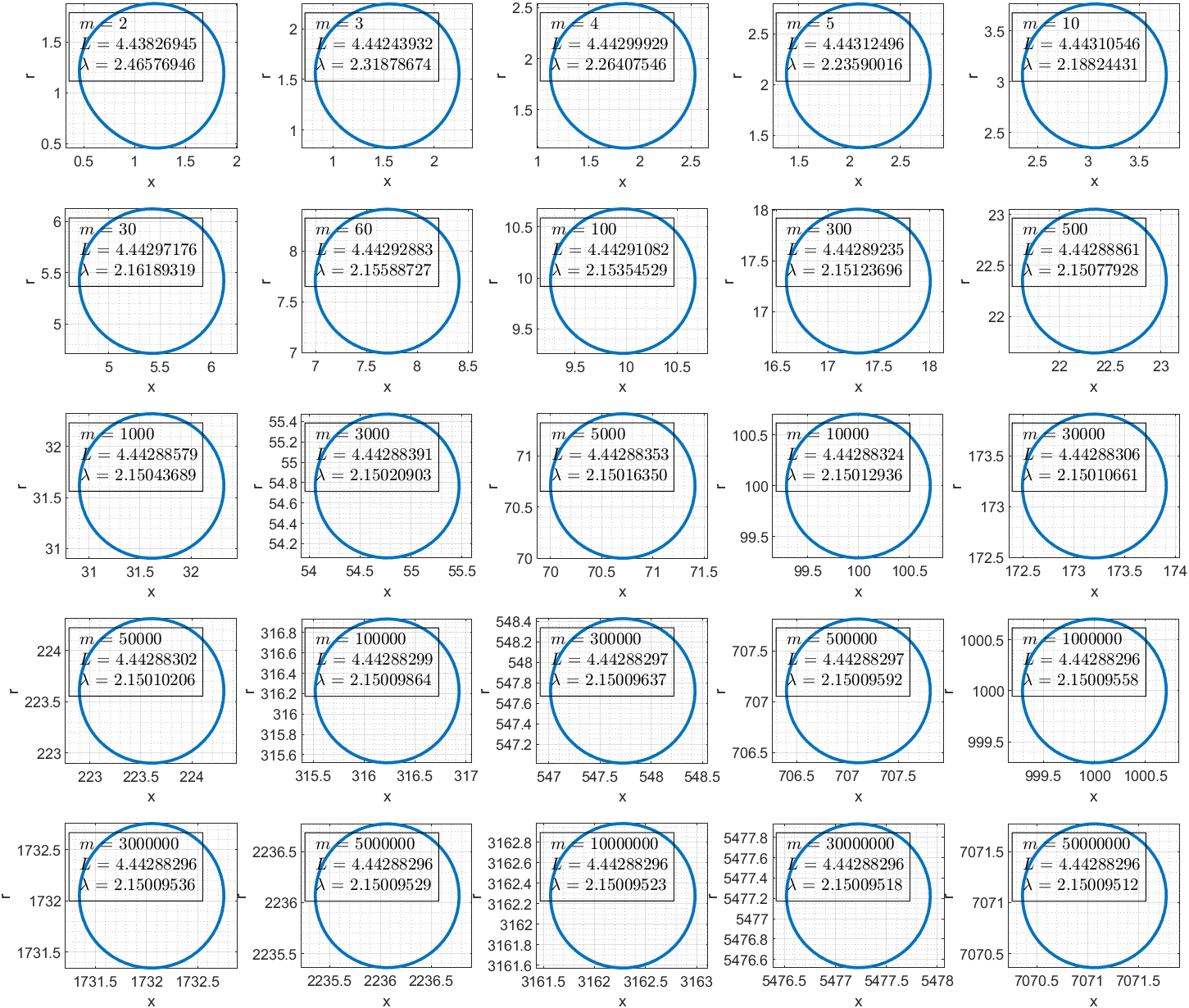}
\caption{ The graph of profile curve of self-shrinkers constructed by Mcgrath. }
\label{fig:Mcgrath_new}
\end{figure}

\begin{figure}[ht]
\includegraphics[width = 0.85\textwidth]{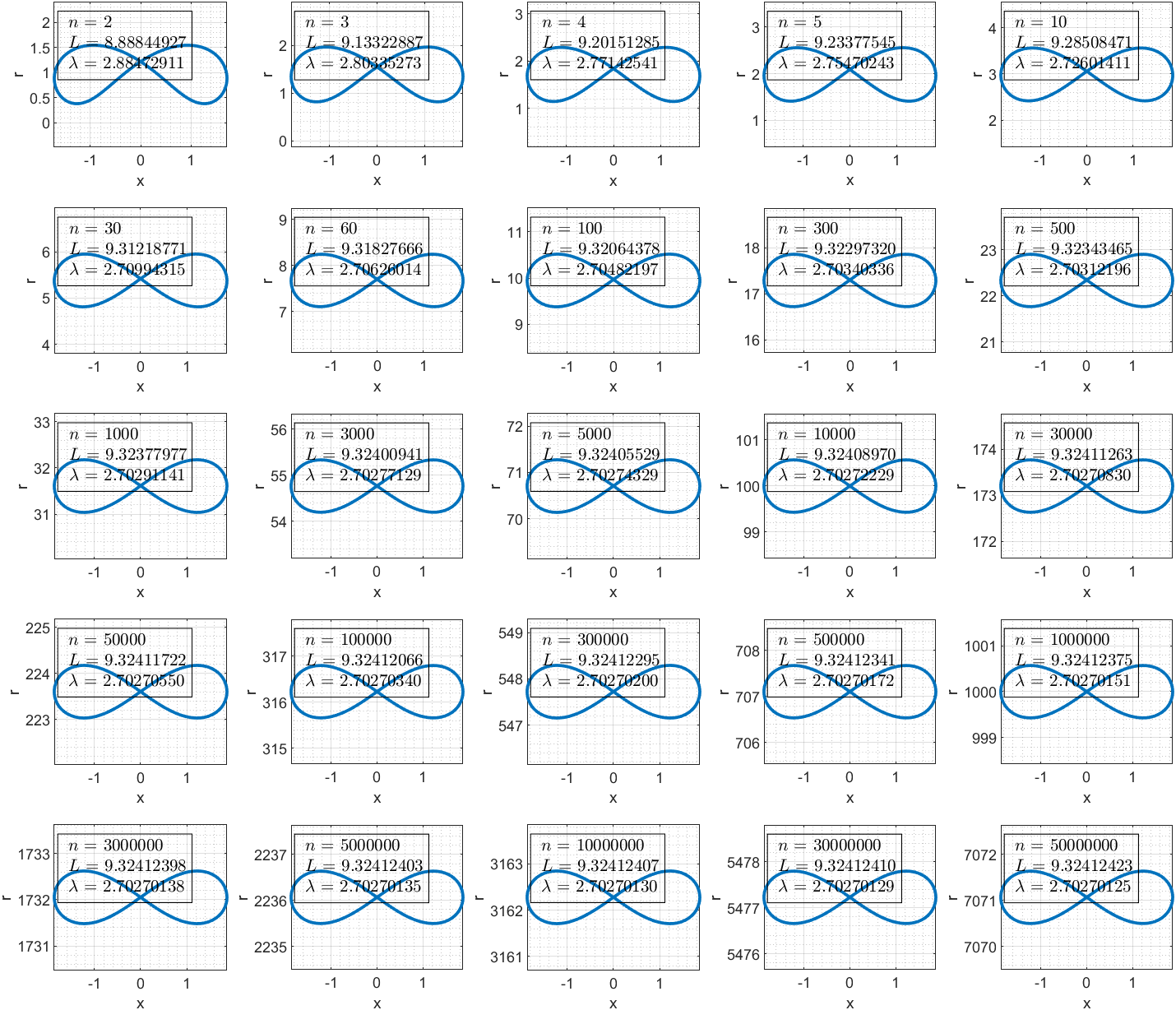}
\caption{ The graph of profile curve of self-shrinkers constructed in \cite{CW1}.}
\label{fig:Angenent_torus_2circles}
\end{figure}

We also numerically estimate the entropy of the self-shrinker which was constructed in \cite{CW1} and obtain that the entropy of 2-dimensional this self-shrinker  is approximately 2.88472911 which is greater than $E_2$. The whole procedure is given in Algorithm \ref{alg:cw1} and the results are shown in Table \ref{tab:Angenent-torus-2circles} and Figure \ref{fig:Angenent_torus_2circles}.

\renewcommand{\algorithmicrequire}{\textbf{ Input:}}
\renewcommand{\algorithmicensure}{\textbf{Output:}}
\begin{algorithm}[ht]
\caption{ Compute the perimeters and entropies for the case in \cite{CW1}}
\label{alg:cw1}
\begin{algorithmic}[1]
\REQUIRE the dimension $n$
\STATE \textbf{Initialize} $a_L = -\frac{\pi}{3}, a_R = 0, L_{max} = 10, tspan = [0,L_{max}]$.
\STATE Set $opts = \mathrm{odeset}(``Events", @eventFun, ``RelTol", 10^{-10})$.
\WHILE{$|a_L-a_R|>10^{-12}$}
\STATE $a_0 =\frac{a_L+a_R}{2}, r_L = \sqrt{n-1}+\min\{5,n\}, r_R = \sqrt{n-1}+10^{-8}$
\WHILE{$|r_L-r_R|>10^{-10}$}
\STATE $r_0 = \frac{r_L+r_R}{2}$
\STATE $\mathbf{y_0} = [0;~ r_0;~ a_0; ~0]$
\STATE Solve \eqref{eq:d-theta} by the high-order Runge-Kutta solver \emph{ode89}, i.e.,
\[
[\textbf{t}, \textbf{x}, \textbf{xe}, \textbf{te}] = \emph{ode89}(@\mathrm{func}, tspan, \mathbf{y_0}, opts),
\]
where
$
\mathrm{func}(\textbf{t}, \textbf{y}) =  \left[ \cos(y_3);\sin(y_3); y_1\sin(y_3)+(\frac{(n-1)}{y_2}-y_2)\cos(y_3); \frac{2^{1-\frac{n}{2}}y_2^{n-1}e^{-\frac{y_1^2+y_2^2}{2}}}{\Gamma(\frac{n}{2})}\right].
$
\IF {$all(\mathbf{x_{:,2}})-r_0<0.2)$}
\STATE $r_R = r_0$;
\ELSE
\IF {$\mathbf{xe}_{end,2}-r_0 > 10^{-10}$}
\STATE $r_L = r_0$
\ELSIF {$r_0 - \mathbf{xe}_{end,2} > 10^{-10}$}
\STATE $r_R = r_0$
\ELSE
\STATE $\mathbf{break}$
\ENDIF
\ENDIF
\ENDWHILE
\STATE $\mathbf{t}_0 = \mathbf{t}, \mathbf{x}_0 = \mathbf{x}, \mathbf{xe}_0 = \mathbf{xe}, \mathbf{te}_0 = \mathbf{te}$
\STATE $\mathbf{y}_0 = (\mathbf{xe}_{end,:})^T, tspan = [\mathbf{t}(end), L_{max}]$
\STATE Solve
$
[\textbf{t}, \textbf{x}, \textbf{xe}, \textbf{te}] = \emph{ode89}(@\mathrm{func}, tspan, \mathbf{y_0}, opts),
$
\IF{$a_0 - \mathrm{mod}(\mathbf{xe}_{end,3}, -\pi) > 10^{-10}$}
\STATE $a_R = a_0$
\ELSIF{$\mathrm{mod}(\mathbf{xe}_{end,3}, -\pi) - a_0> 10^{-10}$}
\STATE $a_L = a_0$
\ELSE
\STATE $\mathbf{break}$
\ENDIF
\ENDWHILE
\STATE Set the perimeter $L = \textbf{t}_{end}$
\STATE Set the entropy $\lambda = \textbf{x}_{end,4}$
\ENSURE the perimeter $L$ and the entropy $\lambda$
\end{algorithmic}
\end{algorithm}

\begin{table}[ht]
    \centering
    \caption{Numerical results for: self-shrinkers constructed in \cite{CW1}}
    \label{tab:Angenent-torus-2circles}
    \resizebox{0.65\linewidth}{!}{
    \begin{tabular}{r|c|c}
   \noalign{\hrule height 1pt}
  Dimension & Perimeter of profile curve  & Entropy  \\
  \noalign{\hrule height 1pt}
2          & 8.88844927 &	2.88472911 \\ \hline
3          & 9.13322887 &	2.80335273 \\ \hline
4          & 9.20151285 &	2.77142541 \\ \hline
5          & 9.23377545 &	2.75470243 \\ \hline
10         & 9.28508471 &	2.72601411 \\ \hline
30         & 9.31218771 &	2.70994315 \\ \hline
60         & 9.31827666 &	2.70626014 \\ \hline
100        & 9.32064378 &	2.70482197 \\ \hline
300        & 9.32297320 &	2.70340336 \\ \hline
500        & 9.32343465 &	2.70312196 \\ \hline
1000       & 9.32377977 &	2.70291141 \\ \hline
3000       & 9.32400941 &	2.70277129 \\ \hline
5000       & 9.32405529 &	2.70274329 \\ \hline
10000      & 9.32408970 &	2.70272229 \\ \hline
30000      & 9.32411263 &	2.70270830 \\ \hline
50000      & 9.32411722 &	2.70270550 \\ \hline
100000     & 9.32412066 &	2.70270340 \\ \hline
300000     & 9.32412295 &	2.70270200 \\ \hline
500000     & 9.32412341 &	2.70270172 \\ \hline
1000000    & 9.32412375 &	2.70270151 \\ \hline
3000000    & 9.32412398 &	2.70270138 \\ \hline
5000000    & 9.32412403 &	2.70270135 \\ \hline
10000000   & 9.32412407 &	2.70270130 \\ \hline
30000000   & 9.32412410 &	2.70270129 \\ \hline
50000000   & 9.32412423 &	2.70270125 \\ 
\noalign{\hrule height 1pt}
\end{tabular}
}
\end{table}

\end{document}